\documentclass[12pt]{amsart}
\usepackage{a4wide}
\usepackage{amssymb,amsmath,amsthm,latexsym}
\usepackage{graphics}
\usepackage{verbatim}
\usepackage{diagrams}

\newcommand{\D}[1]{{\it #1}}

\newtheorem{theorem}{Theorem}[section]
\newtheorem{lemma}[theorem]{Lemma}
\theoremstyle{remark}
\newtheorem{remark}[theorem]{Remark}
\newtheorem{example}[theorem]{Example}

\newcommand{\set}[1]{\ensuremath{\left\{#1 \right\}}}
\newcommand{\setcond}[2]{\set{#1 \; \left| \: #2 \right.}}
\newcommand{\linhull}[1]{\ensuremath{\left\langle#1 \right\rangle}}
\newcommand{\lhcond}[2]{\linhull{#1 \; \left| \: #2 \right.}}
\newcommand{\Sum}{\sum\limits}
\newcommand{\Prod}{\prod\limits}

\newcommand{\scxf}[1]{\mathbf{#1}}
\newcommand{\pruh}{\overline}
\newcommand{\full}[1]{\scxf 2^{#1}}

\newcommand{\cpxf}[1]{\mathcal{#1}}
\newcommand{\pd}{\partial}
\newcommand{\difr}{\hbox{d}}
\newcommand{\CxrH}{\cpxf{\tilde C}_{\circledast}}
\newcommand{\CxrC}{\cpxf{\tilde C}^{\circledast}}

\newcommand{\Hred}{\tilde H} 
\newcommand{\HrSC}[2]{\Hred_{#1} \left(\scxf{#2}\right)} 
\newcommand{\HrFree}[2]{\Hred_{#1} \left({#2}\right)}
\newcommand{\HrPair}[3]{\Hred_{#1} \left(#2, #3 \right)}
\newcommand{\norm}[1]{\left|\kern-.3mm\left| #1 \right|\kern-.3mm\right| }

\newcommand{\rarw}{\rightarrow}

\newcommand{\im}{\mathop{\mathrm{im}}}

\newcommand{\sgn}{\mathop{\mathrm{sgn}}}
\newcommand{\zet}{\mathbb{Z}} 
\newcommand{\ring}{R}

\newcommand{\twolines}[2]{\genfrac{}{}{0pt}{}{#1}{#2}}
\newcommand{\isomorphic}{\cong}

\newcommand{\mybibitem}[5]{\bibitem{#1}{#2: #3, \textsl{#4} #5.}}

\begin{document}
\title[Combinatorial Alexander Duality]
{Combinatorial Alexander Duality \\--- A Short and Elementary Proof}
\author{Anders Bj\"orner 
and Martin Tancer}

\address{Department of Mathematics, Royal Institute of Technology,
S-100 44 Stockholm, Sweden}
\email{bjorner@math.kth.se}

\address{Department of Applied Mathematics, Faculty of Mathematics and Physics, Charles
University, Ma\-lo\-stran\-sk\'e n\'am\v{e}st\'i~25, 118~00, Prague, Czech Republic.}
\email{martin@atrey.karlin.mff.cuni.cz}

\maketitle

\begin{abstract}
Let $\scxf X$ be a simplicial complex with ground set $V$. Define its
Alexander dual as the simplicial complex $\scxf X^* = \setcond{\sigma \subseteq V}{V
\setminus \sigma \notin \scxf X}$. The combinatorial Alexander duality states
that the $i$-th reduced homology group of $\scxf X$ is isomorphic to the
$(|V|-i-3)$-th reduced cohomology group of  $\scxf X^*$ (over a given
commutative ring $\ring$). We give a self-contained proof from first principles,
accessible to the non-expert.
\end{abstract}

\section{Introduction}

Let $\scxf X$ be a simplicial complex with ground set $V$. For $\sigma \in \scxf X$
let $\pruh\sigma = V \setminus \sigma$. The \D{Alexander dual} of $\scxf X$ is 
the simplicial complex on the same ground set defined by 
$$
\scxf X^* = \left(V, \setcond{\sigma \subseteq V}{\pruh\sigma \notin \scxf
X} \right).
$$
See Figure~\ref{figS} for an example of a simplicial complex and its Alexander
dual.  

\begin{figure}[h]
\centering
\includegraphics{dual.1}
\caption{Simplicial complex $\scxf S$ and its dual.}
\label{figS}
\end{figure}

It is easy to see that $\scxf X^{**} = \scxf X$. Furthermore, a close homological connection 
exists between $\scxf X$ and   $\scxf X^{*}$, that in the combinatorics folklore has become known
as  ``combinatorial Alexander duality''.
It states that knowledge of the homology of a simplicial
complex gives knowledge of the cohomology of its Alexander dual:

\begin{theorem}[Combinatorial Alexander Duality]
\label{ADThm}
Let $\scxf X$ be a simplicial complex with a ground set of the size $n$. Then
$$
\HrSC i{X} \isomorphic \Hred^{n-i-3} (\scxf X^*).
$$
(Here $\Hred$
stands for reduced homology resp. cohomology over a given ring $\ring$.)
\end{theorem}

The earliest explicit statements of Theorem \ref{ADThm} that we know of appear in
Kalai \cite[p. 348]{Kal} and Stanley \cite[p. 184]{Sta3}.
Combinatorial Alexander duality is a special case of the original Alexander
duality:

\begin{theorem} \rm{(Alexander \cite{Ale}, 1922)}
 Let $A$ be a subset of the sphere $S^n$, such that the pair $(S^n, A)$
 is triangulable. Then
 $$
 \HrFree iA \isomorphic \Hred^{n-i-1}(S^n \setminus A).
 $$
\end{theorem}

The connection is the following: Suppose that $\scxf X$ is a simplicial complex
different from the full simplex with ground set $\set{1, 2, \dots, n + 2}$. Let
$\scxf Y \isomorphic S^n$ be the $n$-skeleton of the full simplex on the set
$\set{1, 2, \dots n + 2}$. Let us denote by $\norm{\scxf X}$ (resp. $\norm{\scxf Y}$) a
geometric realization of $\scxf X$ (resp. $\scxf Y$). Then $\norm{\scxf X}
\subseteq \norm{\scxf Y}$, and it can be shown that 
$\norm{\scxf Y} \setminus \norm{\scxf X}$ is 
homotopy equivalent to $\norm{\scxf X^*}$. Thus if $A = \norm{\scxf X}$ is seen as a
subset of $S^n$ we get 
$$\HrSC i{X} \isomorphic 
\HrFree iA \isomorphic \Hred^{n-i-1} (S^n \setminus A)  \isomorphic \Hred^{(n+2)-i-3} (\scxf X^*).
$$

The Alexander duality theorem has played a very important role in the development
of algebraic topology. See \cite{Hat} or \cite{Mun} for 
context and modern treatments, and \cite{Jam}
for interesting historical information about James Waddell Alexander, the man and his mathematics.
\medskip

The modest task of this article is to make the ideas behind combinatorial Alexander duality
more widely accessible by giving a ``combinatorial'' proof from first principles, not relying on 
more general techniques.
Central to our approach is a poset point of view and scrutiny of the combinatorics of the sign labeling
of edges in the Boolean lattice $\full V$ used for the (co)boundary operations in
(co)homology.
Applications of combinatorial Alexander duality  in combinatorics and algebra can be found for
example in \cite{BBM}, \cite{BZ}, \cite{BP}, \cite{ER}, 
\cite{Kal}, \cite{KalMes}, \cite{MilStu}, \cite{Sta3}, \cite{Sta4}, \cite{Sta1}, \cite{Sta2}.

Combinatorial Alexander duality exists in more general versions. 
One such generalization involves induced subcomplexes of $\scxf X$
versus links of faces of $\scxf X^*$, see \cite[p. 28]{BP}. For another one,
let $\scxf X$ be a subcomplex of the boundary complex of a 
$(n-1)$-dimensional convex polytope $P$
and let $\scxf X^*$  be the subcomplex of the boundary of the
dual polytope $P^*$ consisting of faces dual to those faces of $P$
that are not in $\scxf X$. Then Theorem \ref{ADThm} 
still holds, see e.g. \cite{BZ} for an application. Furthermore, a very general 
Alexander duality theorem for nonacyclic
Gorenstein complexes appears in \cite[p. 66]{Sta1}. The reason that we confine this article
to the particular simplicial setting of Theorem \ref{ADThm}, although the idea of
the proof  is correct in greater generality, is that we want to work with the explicit
sign labeling offered by the environment of the Boolean lattice.

\section{Preliminaries}
\label{sechc}

Let $\scxf X$ be a simplicial complex with ground set $V = \set{1,2, \dots,
n}$. For $j \in \sigma \in \scxf X$ we define the {\em sign} $\sgn(j, \sigma)$ as
$(-1)^{i-1}$,
where $j$ is the $i$-th smallest element of the set $\sigma$. 
The following simple property of the sign function will be needed.
\begin{lemma}
\label{signs}
Let $k\in \sigma  \subseteq \set{1, 2, \dots, n}$ and  $
p( \sigma) = \Prod_{i \in \sigma}(-1)^{i-1}.
$
Then
$$
\sgn(k, \sigma)\, p(\sigma\setminus k) = \sgn (k, \pruh\sigma \cup k)\, p(\sigma).
$$
\end{lemma}
\begin{proof}
We have that
$$\sgn (k, \sigma)\, \sgn (k, \pruh\sigma \cup k)=
\Prod_{\twolines{i \in \sigma}{i < k}}(-1) \,
\Prod_{\twolines{i \in \pruh\sigma}{i<k}}(-1)=(-1)^{k-1}
$$
and
$$ p(\sigma)\, p(\sigma\setminus k) =
\Prod_{i \in \sigma}(-1)^{i-1} \,
\Prod_{i \in \sigma\setminus k}(-1)^{i-1} =(-1)^{k-1}
$$
\end{proof}

In the rest of this section we review the definitions and notation used for 
(co)homology.
Throughout the paper suppose that $R$ is a commutative ring containing 
a unit element. 

\subsection{Reduced Homology}
\label{redhom}
 Let $C_i = C_i (\scxf X)$ be a free
$\ring$-module with the free basis $\setcond{e_\sigma}{\sigma \in \scxf X,
\dim \sigma = i}$. The \D{reduced chain complex} of $\scxf X$ over $\ring$ is
the complex 
$$  
   \CxrH(\scxf X) = \CxrH ({\scxf X}; \ring) = \cdots 
   \lTo C_{i-1} \lTo^{\textstyle\pd_i}
   C_i \lTo^{\textstyle\pd_{i+1}} 
   C_{i+1} \lTo \cdots, \qquad i \in \zet,
$$
whose mappings $\pd_i$ are defined as
$$
\pd_i(e_\sigma) = \Sum_{j \in \sigma} \sgn (j, \sigma) e_{\sigma \setminus j}.
$$
The complex
$\CxrH (\scxf X)$ is formally infinite; however, $C_i = 0$ for $i <
-1$ or~$i > \dim \scxf X$.
The \D{$n$-th reduced homology group} of $\scxf X$ over $\ring$ is defined as
$$
\HrSC nX = \Hred_n \left( \scxf X; \ring \right) = \ker \pd_n / \im \pd_{n+1}
.
$$

\subsection{Reduced Cohomology}
\label{redcohom}
Let $C^i = C^i(\scxf X)$ be a free $\ring$-module with the free basis
$\setcond{e^*_\sigma}{\sigma \in \scxf X, \dim \sigma = i}$. The \D{reduced
cochain complex} of $\scxf X$ over $\ring$ is the complex
$$
\CxrC(\scxf X) =  \CxrC ({\scxf X}; \ring) = \cdots 
   \rTo  C^{i-1} \rTo^{\textstyle\pd^i} C^i
   \rTo^{\textstyle\pd^{i+1}} C^{i+1} \rTo \cdots, 
   \qquad i \in \zet,
$$
where $\pd^i = \pd^*_i$ are maps dual to $\pd_i$, explicitly stated:
$$
\pd^i(e^*_\sigma) = \Sum_{\twolines{j \notin \sigma}{\sigma \cup j \in \scxf X}}
\sgn (j, \sigma \cup j) e^*_{\sigma \cup j}.
$$
The \D{$n$-th reduced cohomology group} of $\scxf X$ over $\ring$ is defined as
$$
\Hred^n \left( \scxf X \right) = \Hred^n \left( \scxf X; \ring \right) = \ker \pd^{n+1} / \im \pd^n.
$$

\subsection{Relative Homology}

Suppose that $\scxf X$ is a simplicial complex and $\scxf A$ is a subcomplex of
$\scxf X$. Let $R_i = R_i( \scxf X, \scxf A) = C_i(\scxf X) / C_i(\scxf A)$,
where $C_i$ was defined in Section~\ref{redhom}.
The \D{relative reduced chain complex} of $(\scxf X, \scxf A)$ over $\ring$ is
the complex
$$
   \CxrH ({\scxf X}, {\scxf A}) =  
   \CxrH ({\scxf X}, {\scxf A}; \ring) = \cdots 
   \lTo R_{i-1} \lTo^{\textstyle\difr_i}
   R_i \lTo^{\textstyle\difr_{i+1}} R_{i+1} 
   \lTo \cdots, \qquad i \in \zet,  
$$
where $\difr_i$ are defined as
$$
\difr_i(e_\sigma + C_i(\scxf A)) = \Sum_{j \in \sigma} \sgn (j, \sigma)\left(
e_{\sigma \setminus j} + C_{i-1}(\scxf A)\right).
$$

The \D{$n$-th relative reduced homology group} of $(\scxf X, \scxf A)$ over
$\ring$ is defined as
$$
\HrFree n{\scxf X, \scxf A} = \Hred_n \left( \scxf X, \scxf A; \ring \right) = \ker \difr_n / \im \difr_{n+1}.
$$

\begin{remark}
\label{relhom}
When we wish to compute relative homology groups, we can identify $R_i = C_i( \scxf X) /
C_i(\scxf A)$ with a free $\ring$-module with the free basis
$\setcond{e_\sigma}{\sigma \in \scxf X, \sigma \notin \scxf A, \dim \sigma = i}$.
Then $\difr_i$ can be rewritten as:
$$
\difr_i(e_\sigma) = \Sum_{\twolines{j \in \sigma}{\sigma \setminus j \notin \scxf
A}} \sgn (j, \sigma)
e_{\sigma \setminus j}.
$$
\end{remark}


One of the important properties of relative homology groups is that they fit
into a~long exact sequence.
See e.~g.~\cite{Hat} or \cite{Mun} for a proof and more details. 
\begin{lemma}[Long Exact Sequence of a Pair]
\label{lespair}
Suppose that $\scxf X$ and $\scxf A$ are simplicial complexes $\scxf A \subseteq
\scxf X$. Then there is a long exact sequence
$$
\cdots \rTo \HrSC nA \rTo \HrSC nX \rTo \HrPair n{\scxf X}{\scxf A}
\rTo \HrSC {n-1}A \rTo \cdots
$$
\end{lemma}

\section{The Idea of the Proof}
\label{secidea}
Before embarking on the proof of Theorem~\ref{ADThm}, we first
present the idea. 

Suppose that $\scxf X$ is a simplicial complex with ground set $V$. Let
$\Gamma$ be the lattice of all subsets of $V$, and let $\Gamma_{\scxf X}$ be the
subposet of $\Gamma$ corresponding to the subsets that are in~$\scxf X$. Then
the $n$-th homology group of $\scxf X$ depends just on the $n$-th and $(n +
1)$-st levels of the poset $\Gamma_{\scxf X}$. Let $\full V$ be a full simplex with
vertex set $V$. It is easy to see that $\HrSC iX \isomorphic \HrPair
{i+1}{\full V}{\scxf X}$ (see Lemma~\ref{ADfirst}). Thus,
we restate the problem as computing homologies of the chain complex determined
by the complement of $\Gamma_{\scxf X}$ --- in the sense of Remark~\ref{relhom}.

The idea of the proof is that if we turn the lattice
upside down (exchange $\sigma \subseteq V$ with its complement), then this
combinatorial map on generators should induce a
canonical isomorphism between the relative homology of the pair $(\full V,
\scxf X)$ and the cohomology of $\scxf X^*$. This idea is basically correct;
however, the mapping is not an isomorphism as described 
--- some sign operations are necessary.

\begin{example}
Let $\scxf S$ be the simplicial complex in Figure~\ref{figS}. Its ground set is
the set $V_{\scxf S} = \set{1, 2, 3, 4}$. The posets $\Gamma_{\scxf S}$ and
$\Gamma_{\scxf S^*}$ are depicted in Figure~\ref{figGammaS}, and the left part
of the picture also shows the complement of $\Gamma_{\scxf S}$ (bold, dashed)
determining the homology of $(\full {V_{\scxf S}}, \scxf S)$.

In the sense of Remark~\ref{relhom}, the chain complex $\CxrH \left( \full
{V_{\scxf S}}, \scxf S \right)$ is 
$$
\cdots \lTo 0 \lTo \linhull{e_{24}, e_{34}}\lTo^{\textstyle\difr_2}
\linhull{e_{123}, e_{124}, e_{134}, e_{234}}
\lTo^{\textstyle\difr_3}
\linhull{e_{1234}} \lTo 0 \lTo \cdots
$$
and the cochain complex $\CxrC \left( \scxf S^* \right)$ is 
$$
\cdots \rTo  0 
   \rTo  \linhull{e^*_{0}}
  \rTo^{\textstyle\pd^{0}} \linhull{ e^*_{1}, e^*_{2}, e^*_{3}, e^*_{4}}
  \rTo^{\textstyle\pd^{1}} \linhull{ e^*_{12}, e^*_{13} }
  \rTo 0 \rTo \cdots
$$  
The map $e_{\sigma} \rarw e^*_{\pruh{\sigma}}$ is {\em not} an isomorphism of these
two complexes (if $2r \neq 0$ for $0 \neq r \in \ring$), since
$\difr_2 (e_{234}) = -e_{24} + e_{34}$, while $\pd^1 (e^*_1) = - e^*_{12} -
e^*_{13}$.  Nevertheless, these two chain complexes are isomorphic, as will be
shown in the next section.
\end{example}

\begin{figure}
\centering
\includegraphics{lattices.1}
\caption{The lattices $\Gamma_{\scxf S}$ (left) and $\Gamma_{\scxf S^*}$ (right).}
\label{figGammaS}
\end{figure}

\section{The Proof}
\label{secprf}

The proof of Theorem~\ref{ADThm} is obtained by combining 
Lemma~\ref{ADfirst} and Lemma~\ref{ADsecond}.

\begin{lemma}
\label{ADfirst}
Let $\scxf X$ be a simplicial complex with ground set $V$. Then
$$
\HrSC iX \isomorphic \HrPair {i+1}{\full V}{\scxf X}.
$$
\end{lemma}

\begin{proof}
This follows from Lemma~\ref{lespair}. There is the long exact sequence of
the pair $(\full V, \scxf X)$:
$$
   \cdots  \rTo  \HrFree {i+1}{\full V}  \rTo  
   \HrPair {i+1}{\full V}{\scxf X}  \rTo  \HrSC iX  \rTo 
   \HrFree i{\full V} \rTo \cdots 
$$
The groups $\HrFree {i+1}{\full V}$ and $\HrFree i{\full V}$ are zero, hence 
 the groups $\HrPair {i+1}{\full V}{\scxf X}$ and~$\HrSC iX$ are isomorphic.
\end{proof}

\begin{lemma}
\label{ADsecond}
Let $\scxf X$ be a simplicial complex with ground set $V$ of size $n$. Then 
$$\HrPair {i+1}{\full V}{\scxf X} \isomorphic \Hred^{n-i-3}(\scxf X^*).$$
\end{lemma}
\begin{proof}

Suppose that $V = \set{1, 2, \dots, n}$.
The chain complex for reduced homology of the pair $(\full V, \scxf X)$ is the
complex
\begin{diagram}
   \qquad & \cdots & \lTo^{\difr_{j-1}} & R_{j-1} & \lTo^{\difr_j} 
   & R_j & \lTo^{\difr_{j+1}} & \cdots, & \qquad & j \in \zet,  \\
\end{diagram}
where 
$
R_j = \lhcond{e_{\sigma}}{\sigma \subseteq V, \sigma \notin \scxf X, \dim
\sigma = j}
$ 
and $\difr_j$ are the unique homomorphisms satisfying
$$
\difr_j (e_\sigma) = 
\Sum_{\twolines{k \in \sigma}{\sigma \setminus k \notin \scxf X}} 
\sgn(k, \sigma) e_{\sigma \setminus k}.
$$
The cochain complex for reduced cohomology of $\scxf X^*$ is the complex
\begin{diagram}
   \qquad & \cdots & \rTo^{\pd^{j-1}} & C^{j-1} & \rTo^{\pd^j} 
   & C^j & \rTo^{\pd^{j+1}} & \cdots, & \qquad & j \in \zet,  \\
\end{diagram}
where
$
C^j = \lhcond{e^*_\sigma}{\sigma \subseteq V, \dim \sigma = j, \sigma
\in \scxf X^*} = 
\lhcond{e^*_\sigma}{\sigma \subseteq V, \dim \pruh\sigma = n-j-2, \pruh\sigma
\notin \scxf X}
$
and $\pd^j$ are the unique homomorphisms satisfying
$$	
\pd^j \left( e^*_\sigma \right) = 
\Sum_{\twolines{k \notin \sigma}{k \cup \sigma \in \scxf X^*}}
\sgn(k, \sigma \cup k) e^*_{\sigma \cup k} =
\Sum_{\twolines{k \in \pruh\sigma}{\pruh\sigma \setminus k \notin \scxf X}}
\sgn(k, \sigma \cup k) e^*_{\pruh{\pruh\sigma \setminus k}}.
$$

Define $p(\sigma)$ as in Lemma \ref{signs}
and let $\phi_j : R_j \rarw C^{n - j - 2}$ be the isomorphisms generated by
the formula
$$\phi_j (e_\sigma) = p(\sigma)\, e^*_{\pruh\sigma}$$
for $\sigma \notin \scxf X$, $\dim \sigma
= j$ (note that these two conditions  are equivalent to $\dim \pruh\sigma = n
- j - 2$, $\pruh\sigma \in \scxf X^*$). We then have a diagram:

\begin{diagram}
\cdots         & \lTo^{\difr_{j-1}}  & R_{j-1}           & \lTo^{\difr_j}     & 
R_j            & \lTo^{\difr_{j+1}}  &     \cdots \\
               &                     & \dTo^{\phi_{j-1}} & \circ              &
 \dTo^{\phi_j} &                     &            \\
\cdots         &\lTo^{\pd^{n - j}}   & C^{n-j-1}         & \lTo^{\pd^{n-j-1}} &
C^{n-j-2}      &\lTo^{\pd^{n-j-2}}   &     \cdots \\
\end{diagram}
\medskip

We check that $\phi_{j-1} \circ \difr_j = \pd^{n-j-1} \circ \phi_j$. Let $\sigma
\subseteq V$, $\sigma \notin \scxf X$, $\dim \sigma = j$. Then
$$
\phi_{j-1} \circ \difr_j (e_{\sigma}) = 
\phi_{j-1} 
\Biggl(\Sum_{\twolines{k \in \sigma}{\sigma \setminus k \notin \scxf
X}} \sgn(k, \sigma) e_{\sigma \setminus k}\Biggr) =
\Sum_{\twolines{k \in \sigma}{\sigma \setminus k \notin \scxf
X}} \sgn(k, \sigma) p(\sigma \setminus k)
e^*_{\pruh{\sigma \setminus k}},
$$
$$
\pd^{n-j-1} \circ \phi_j (e_{\sigma}) =
\pd^{n-j-1} \left(p(\sigma) \,
e^*_{\pruh \sigma}\right) =
\Sum_{\twolines{k \in
\sigma}{\sigma \setminus k \notin \scxf X}}
p(\sigma)\sgn(k, \pruh\sigma \cup k) 
e^*_{\pruh{\sigma \setminus k}}.
$$
These two sums are equal term by term, due to Lemma~\ref{signs}.
Thus $\phi$ is an isomorphism of the complexes, implying
$$
\HrPair {i+1}{\full V}{\scxf X} \isomorphic 
\Hred^{n-i-3} (\scxf X^*). 
$$
\end{proof}

\noindent
{\sc Acknowledgment:} We are grateful to Ji\v{r}\'{\i}  Matou\v{s}ek and Kathrin Vorwerk 
for helpful remarks.

\end{document}